\documentclass[12 pt]{amsart}
\usepackage{amscd,amsfonts,amssymb,amsmath}
\usepackage{hyperref}
\usepackage{epsfig} 
\usepackage{mathtools}
\usepackage{tabu}
\usepackage{tikz-cd}
\newtheorem{theorem}{Theorem}[section]
\newtheorem{corollary}[theorem]{Corollary}
\newtheorem{lemma}[theorem]{Lemma}
\newtheorem{proposition}[theorem]{Proposition}
\theoremstyle{definition}

\newtheorem{example}[theorem]{Example}
\newtheorem{remark}[theorem]{Remark}

\usepackage[all,cmtip]{xy}
\usepackage{float}
\usepackage{graphicx} 

\begin{document}
\title[Multi-virtual  braid groups]{Multi-virtual  braid groups}
\author[V. Bardakov]{Valeriy G. Bardakov}
\author[T. Kozlovskaya]{Tatyana A. Kozlovskaya}
\author[K. Negi]{Komal Negi}
\author[M. Prabhakar]{Madeti Prabhakar}

\address{Sobolev Institute of Mathematics, 4 Acad. Koptyug avenue, 630090, Novosibirsk, Russia.}
\address{Novosibirsk State Agrarian University, Dobrolyubova street, 160, Novosibirsk, 630039, Russia.}
\address{Regional Scientific and Educational Mathematical Center of Tomsk State University,~36 Lenin Ave., Tomsk, Russia.}
\email{bardakov@math.nsc.ru}

\address{Regional Scientific and Educational Mathematical Center of Tomsk State University,~36 Lenin Ave., Tomsk, Russia.}
\email{t.kozlovskaya@math.tsu.ru}

\address{Indian Statistical Institute, Delhi Center, New Delhi-110016, India}
\email{komal@isid.ac.in}

\address{Indian Institute of Technology, Ropar-140001, Punjab, India}
\email{prabhakar@iitrpr.ac.in}


\begin{abstract}
L. Kauffman (2024) introduced  multi-virtual and  symmetric multi-virtual
 braid groups, which are  generalizations of the  virtual braid group.    We introduce  
 multi-virtual pure and multi-virtual semi-pure braid groups, which are normal subgroups of index $n!$. We give a set of generators and defining relations for these groups,  show that multi-virtual (symmetric multi-virtual) braid group is a semi-direct products of  multi-virtual pure (symmetric multi-virtual pure) braid group and symmetric group.
Also, we introduce  multi-welded and multi-unrestricted braid groups  and some its subgroups and quotients. The paper concludes by outlining open problems and suggesting avenues for future research in this area.  

\end{abstract}

\maketitle 


\section{Introduction}


The braid groups, $B_n$, $n \geq 1$ form a basis of classical knot theory. During some last decades were introduced and intensively studied different generalizations of   $B_n$. For example, virtual braid group $VB_n$, welded  braid group $WB_n$,
singular  braid group $SB_n$ and some other  (see  \cite{B, Bir, KL}  and references therein). These groups form basis of the virtual knot theory,  
 the welded knot theory and the singular knot theory respectively. Also, it is possible to construct some quotient of $VB_n$. On this way,  were defined flat virtual braid group 
$FVB_n$, unrestricted virtual braid group $UVB_n$~\cite{BBD} and some others.

All these groups  $VB_n$,  $WB_n$,   $SB_n$, $FVB_n$, and $UVB_n$ have two types of generators that corresponds to  two types of crossings in the corresponding knot theories. It is natural to call these groups  by 2-sort braid-like groups.
There are braid-like groups, which have three types of generators. For example, singular virtual braid groups~\cite{BK,BK-1}, twisted virtual braid group~\cite{BKNP}. We shall call them by 3-sort braid-like groups. Also, singular twisted virtual braid groups are the 4-sort braid-like groups studied in~\cite{KP}.

L. Kauffman \cite{Ka1} introduced multi-virtual braid groups.  
These groups are  a generalization of the virtual braid group and contains many types of virtual generators. 
Multi-virtual braid groups correspond to multi-virtual knot theory with many types of virtual
crossings. Defining relations in the multi-virtual braid groups are defined by detour moves which  allowed  with respect to some types of crossings, but are forbidden between other crossings. 

The braid group $B_n$  has an epimorphism onto the symmetric group $S_n$, under which the generator $\sigma_i$ goes  to the transposition $(i, i+1)$. The kernel of this epimorphism is said to be a pure braid group.
By analogy, the virtual pure braid group $VP_n$ is the kernel of the epimorphism of the virtual braid group $VB_n$ onto $S_n$ that maps, for each $i$,  $\sigma_i$ and $\rho_i$ to the  $\rho_i$. Here  $\langle \rho_1, \rho_2, \ldots, \rho_{n-1} \rangle \cong S_n$.

In this paper we provide generators and defining relations for pure and semi-pure multi-virtual braid groups and we define some new quotients of these groups, showing generators and defining relations.

The paper is organized as follows. 
In Section \ref{BD}, we recall some known definitions and facts in (virtual) braid
theory. In particular, we give definitions and presentations for virtual pure and virtual semi-pure braid groups. 
In Subsection \ref{MVKT}, we recall multi-virtual knot theory in brief.

In Section \ref{pure}, we recall the Kauffman's definitions  of multi-virtual braid group and symmetric multi-virtual braid group. Further, we introduce multi-virtual pure braid group and multi-virtual semi-pure braid group. We give sets of generators and defining relations for these groups.
Similar results we give for the symmetric multi-virtual braid group.
In Section \ref{QG}, we introduce  multi-welded braid group and multi-unrestricted braid groups, we find their set of generators and defining relation. Also, we defined pure and semi-pure subgroups in these groups.

In the end of the paper, we formulate some open problems and suggest directions for further research.

\section{Basic definitions} \label{BD}

In this section we recall some known definitions which can be found in \cite{Artin, Bir1, Mar}. 

\subsection{Braid group} \label{BG}
\begin{figure}[h]
\includegraphics[totalheight=3cm]{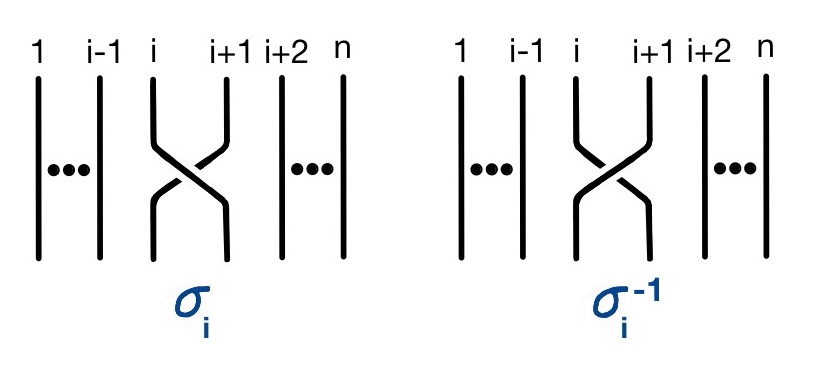}
\caption{The geometric interpretation of $\sigma_i$ and  $\sigma_i^{-1}$} \label{figure1}
\end{figure}

The braid group $B_n$, $n\geq 2$, on $n$ strands can be defined as
a group generated by $\sigma_1,\sigma_2,\ldots,\sigma_{n-1}$ with the defining relations
\begin{equation}
\sigma_i \, \sigma_{i+1} \, \sigma_i = \sigma_{i+1} \, \sigma_i \, \sigma_{i+1},~~~ i=1,2,\ldots,n-2, \label{eq1}
\end{equation}
\begin{equation}
\sigma_i \, \sigma_j = \sigma_j \, \sigma_i,~~~|i-j|\geq 2. \label{eq2}
\end{equation}
The geometric interpretation of  $\sigma_i$ and  its inverse $\sigma_{i}^{-1}$ are depicted  in the Figure~\ref{figure1}.

\subsection{Virtual braid group}  \label{VBG}
The {\it virtual braid group} $VB_n$  was introduced in \cite{Ka}. This group   is generated by the  braid group $B_n$ and the symmetric group $S_n=\langle \rho_1, \rho_2,\ldots, \rho_{n-1} \rangle$  with the following relations:
\begin{align*}
\sigma_i \sigma_{i+1} \sigma_i&=\sigma_{i+1} \sigma_i \sigma_{i+1} & i=1, 2, \ldots, {n-2}, \\
\sigma_i \sigma_j&=\sigma_j \sigma_i &  |i-j| \geq 2, \\
\rho_i^{2}&=e &  i=1, 2, \ldots, {n-1},\\
\rho_i \rho_j&= \rho_j \rho_i &  |i-j| \geq 2,\\
\rho_i \rho_{i+1} \rho_i&= \rho_{i+1} \rho_{i} \rho_{i+1} & i=1, 2, \ldots, {n-2},\\
\sigma_i \rho_j&= \rho_j \sigma_i & |i-j| \geq 2 ,\\
\rho_i \rho_{i+1} \sigma_i&= \sigma_{i+1} \rho_i \rho_{i+1} & i=1, 2, \ldots, {n-2}.
\end{align*}

The virtual pure braid group $VP_n$, $n\geq 2$, was introduced in  \cite{B} as the kernel of the homomorphism $VB_n \to S_n$, $\sigma_i \mapsto \rho_i$, 
$\rho_i \mapsto \rho_i$ for all $i=1, 2, \ldots, i-1$. It was proved that $VP_n$
 admits a
presentation with the  generators $\lambda_{ij},\ 1\leq i\neq j\leq n,$
and the following relations:
\begin{align}
& \lambda_{ij}\lambda_{kl}=\lambda_{kl}\lambda_{ij} \label{rel},\\
&
\lambda_{ik}\lambda_{jk}\lambda_{ij}=\lambda_{ij}\lambda_{jk}\lambda_{ik}
\label{relation},
\end{align}
where distinct letters stand for distinct indices.
The generators of $VP_n$ can be expressed  in terms of  the generators of $VB_n$ by the formulas 
\begin{align*}
\lambda_{i,i+1} &= \rho_i \, \sigma_i^{-1},\\
  \lambda_{i+1,i} &= \rho_i \, \lambda_{i,i+1} \, \rho_i = \sigma_i^{-1} \, \rho_i
\end{align*}
for  $i=1, 2, \ldots, n-1$, and
\begin{align*}
\lambda_{i,j} & = \rho_{j-1} \, \rho_{j-2} \ldots \rho_{i+1} \, \lambda_{i,i+1} \, \rho_{i+1} \ldots \rho_{j-2} \, \rho_{j-1}, \\ 
\lambda_{j,i} & =  \rho_{j-1} \, \rho_{j-2} \ldots \rho_{i+1} \, \lambda_{i+1,i} \, \rho_{i+1} \ldots \rho_{j-2} \, \rho_{j-1}
\end{align*} 
for $1 \leq i < j-1 \leq n-1$.


We have decomposition $VB_n = VP_n \rtimes S_n$ and $S_n$ acts on $VP_n$ by the rules

\begin{lemma}[\cite{B}] \label{form}
Let $a$ be an element of $\langle \rho_1, \rho_2, \ldots, \rho_{n-1} \rangle$ and $\bar{a}$ is its image in $S_n$ under the isomorphism $\rho_i \mapsto (i,i+1)$, $i = 1, 2, \ldots, n-1$, then for any generator $\lambda_{ij}$ of $VP_n$ the following holds
$$
a^{-1} \lambda_{ij} a = \lambda_{(i)\bar{a}, (j)\bar{a}},
$$
where $(k)\bar{a}$ is the image of $k$ under the action of the permutation $\bar{a}$.
\end{lemma}

The flat virtual braid group $FVB_n$ is the quotient of $VP_n$ by the relations $\sigma_i^2 = e$ for all $i$. It contains 
 the group of flat virtual pure braids $FVP_n$ that is generated by elements $\lambda_{ij},\ 1\leq i < j\leq n,$ and  the map 
$\lambda_{ji} \mapsto \lambda_{ij}^{-1}$, $1\leq i < j\leq n,$ defines a homomorphism of $VP_n$ onto $FVP_n$. See Figure~\ref{pure}
\begin{figure}[ht]
\includegraphics[totalheight=7cm]{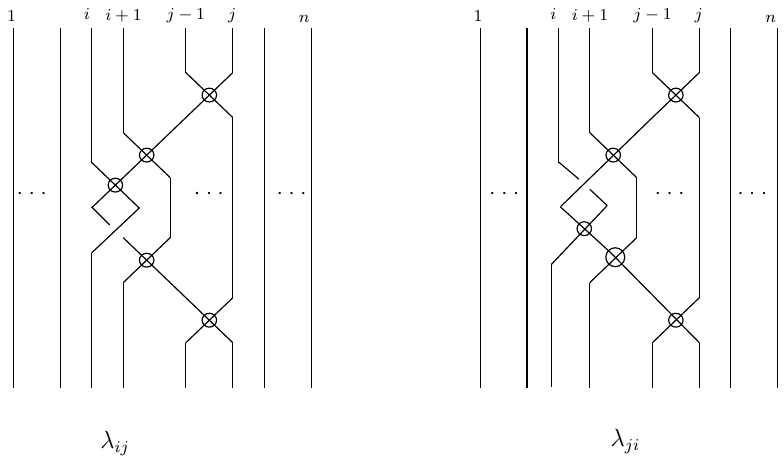}
\caption{Geometric interpretations of generators $\lambda_{ij}$ and $\lambda_{ji}$} \label{pure}
\end{figure}
\begin{example} $VP_3$ is generated by elements
$$
\lambda_{12},~~\lambda_{21},~~\lambda_{13},~~\lambda_{23},~~\lambda_{31},~~\lambda_{32},
$$
and is defined by six relations
$$
\lambda_{ki}\lambda_{kj}\lambda_{ij}=\lambda_{ij}\lambda_{kj}\lambda_{ki},~~i, j, k \in \{1, 2, 3 \}.
$$
The group
$$
FVP_3 = \langle \lambda_{12}, \lambda_{13}, \lambda_{23}~|~  \lambda_{12} \lambda_{13} \lambda_{23} = \lambda_{23} \lambda_{13}  \lambda_{12} \rangle
$$
 with one defining relation, is a subgroup of group $VP_3$. Since, there is a natural retraction $f:VP_3 \to FVP_3$ such that $f\!\mid_{FVP_3}$ is the identity map.
 
 This group is the flat virtual pure braid group and is the  quotient of $VP_3$ by relations $\lambda_{ij} = \lambda_{ji}^{-1}$.
\end{example}

\medskip

It is possible to define another epimorphism $\psi  \colon VB_n \to S_n$ as follows:
$$
\psi (\sigma_i)=e, \; \psi (\rho_i)=\rho_i, \;i=1,2,\dots, n-1.
$$
Let $VH_n$ denote the normal closure of $B_n$ in $VB_n$.
It is evident that $\ker \psi $ coincides with $VH_n$.
Let us define elements:
$$
x_{i,i+1}=\sigma_i,~~x_{i+1,i}=\rho_i \sigma_i \rho_i =\rho_i x_{i,i+1} \rho_i,
$$
for $i= 1, 2, \ldots, n-1$, and 
$$
x_{i,j}=\rho_{j-1} \cdots \rho_{i+1} \sigma_i \rho_{i+1} \cdots \rho_{j-1},
$$
$$
x_{j,i}=\rho_{j-1} \cdots \rho_{i+1} \rho_i \sigma_i \rho_i \rho_{i+1} \cdots \rho_{j-1},
$$
for $1 \le i < j-1 \le n-1$.

The group $VH_n$ admits a presentation with the generators $x_{k,\, l},$ $1 \leq k \neq l \leq
n$,
and the defining relations:
\begin{equation} \label{eq40}
x_{i,j} \,  x_{k,\, l} = x_{k,\, l}  \, x_{i,j},
\end{equation}
\begin{equation} \label{eq41}
x_{i,k} \,  x_{k,j} \,  x_{i,k} =  x_{k,j} \,  x_{i,k} \, x_{k,j},
\end{equation}
where  distinct letters stand for distinct indices.
This presentation was found in \cite{R} (see also \cite{BB}). We have decomposition
$VB_n = VH_n \rtimes S_n$, where  $S_n=\langle \rho_1, \dots, \rho_{n-1} \rangle$ acts on the set
$\{x_{i,j} \, , \;  1 \le i \not= j \le n \}$ by permutation of indices and we have the following  analogous of Lemma \ref{form}.

\begin{lemma}[\cite{BB}] \label{form1}
Let $a$ be an element of $\langle \rho_1, \rho_2, \ldots, \rho_{n-1} \rangle$ and $\bar{a}$ is its image in the symmetric group $S_n$ under the isomorphism $\rho_i \mapsto (i,i+1)$, $i = 1, 2, \ldots, n-1$, then for any generator $x_{ij}$ of $VH_n$ the following holds
$$
a^{-1} x_{ij} a = x_{(i)\bar{a}, (j)\bar{a}},
$$
where $(k)\bar{a}$ is the image of $k$ under the action of the permutation $\bar{a}$.
\end{lemma}


\subsection{Knot theory and multi-virtual knot theory} \label{MVKT}
\begin{figure}[ht]
\includegraphics[totalheight=11cm]{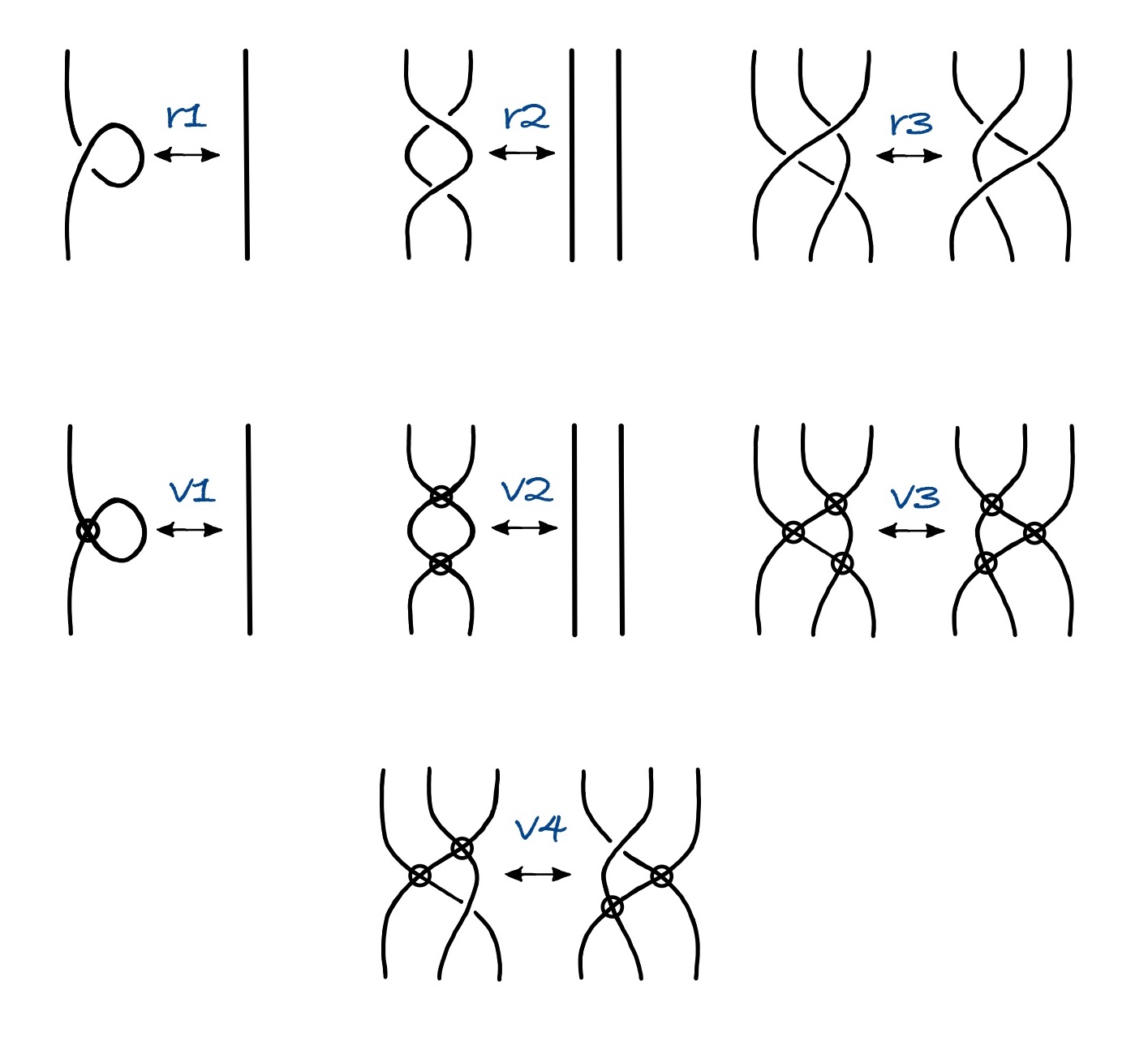}
\caption{Reidemeister and virtual Reidemeister moves} \label{figure4}
\end{figure}

The famous  Reidemeister's theorem~\cite{Reid} says that the set of local  moves  $(r1) - (r3)$ of link's diagrams, shown in the first row in Figure~\ref{figure4} define equivalent relation on link's diagrams. Reidemeister showed that if we have two knots (links) in 3-dimensional space, then they are ambient isotopic if and only if corresponding diagrams for them can be obtained, one from the other, by a finite sequence of moves $(r1) - (r3)$. 

A virtual link represents a natural combinatorial generalization of a classical link (see~\cite{Ka}).
A virtual diagram (or a diagram of a virtual link) is the image of an immersion of a framed 4-valent graph in $\mathbb{R}^2$ with a finite number of intersections of edges. Moreover, each intersection is a transverse double point which we call a virtual crossing and mark by a small circle, and each vertex of the graph is endowed with the classical crossing structure (with a choice for underpass and overpass specified). The vertices of the graph are called classical crossings or just crossings.

\begin{figure}[h]
\noindent\centering{\includegraphics[height=0.3 \textwidth]{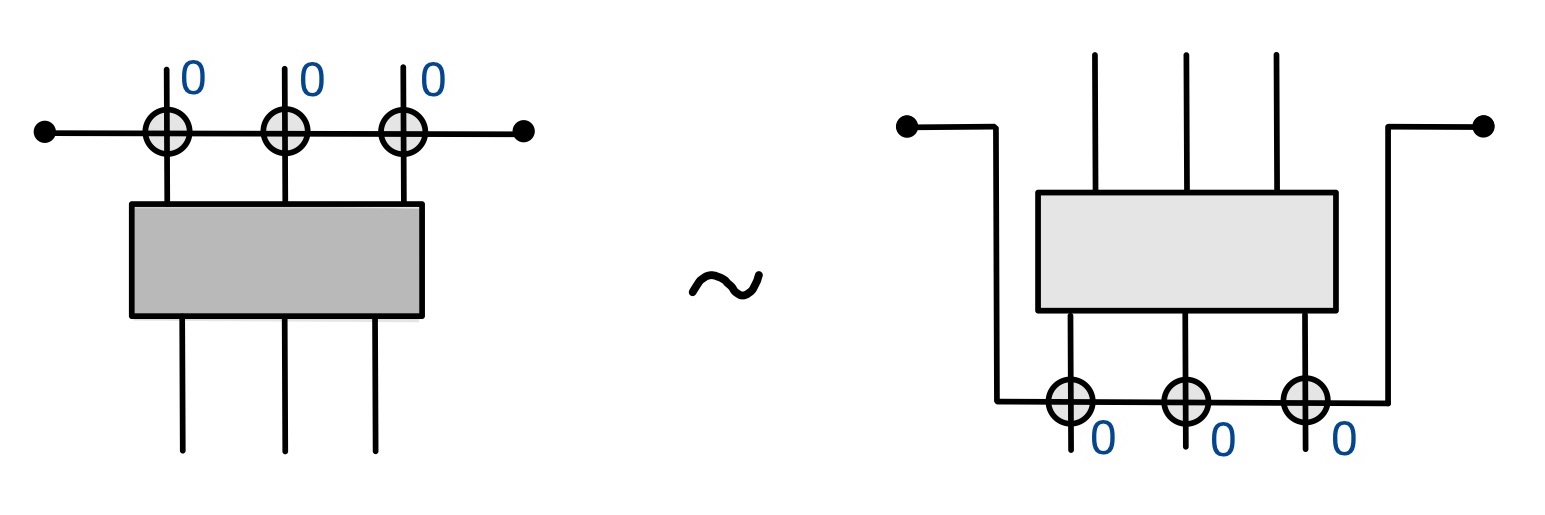}}
\caption{Detour move for a virtual crossings of type $0$} \label{det}
\end{figure}

Two virtual link diagrams are virtually equivalent if one diagram can be transformed into the other by a sequence of classical and virtual Reidemeister moves as shown in Figure ~\ref{figure4}. The virtual Reidemeister moves are equivalent to a single move, the detour move, which is executed by selecting a segment of a component of the link diagram that contains no classical crossings (see Figure~\ref{det}). After removing this segment, we may insert a new segment with no triple points and any double points result in a new virtual crossing.
A virtual link is an equivalence class of virtual diagrams modulo the  Reidemeister moves referring to classical crossings and the detour moves.

\begin{figure}[h]
\noindent\centering{\includegraphics[height=0.3 \textwidth]{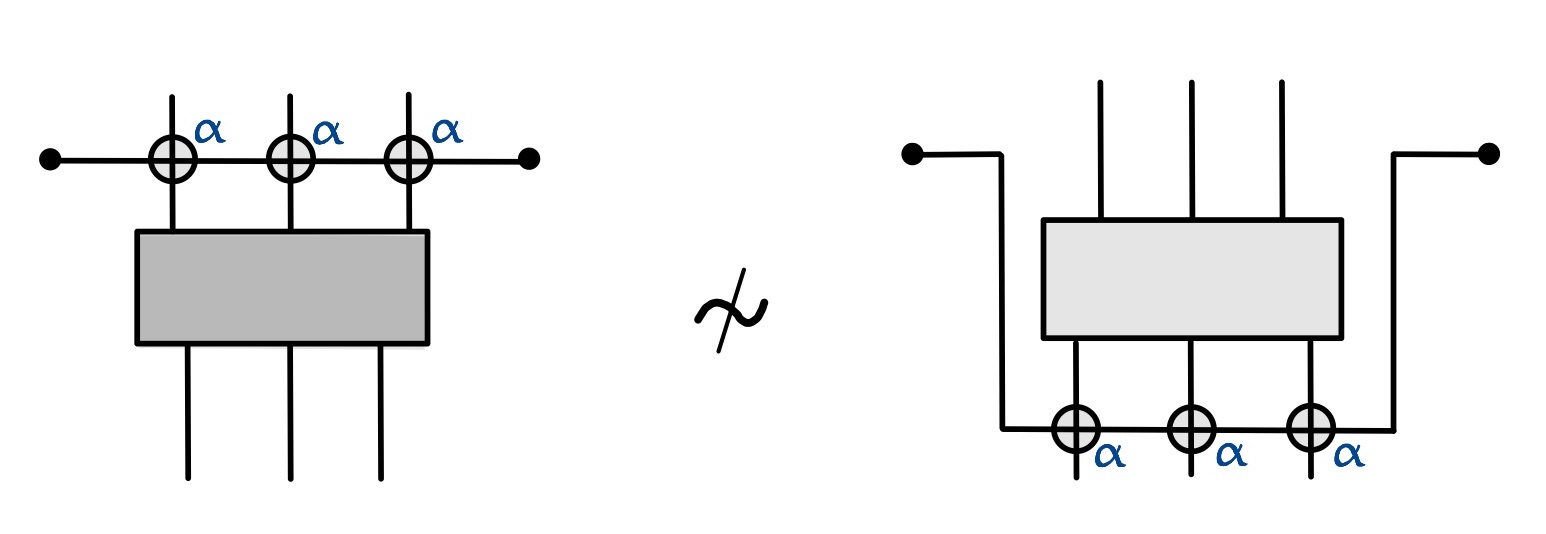}}
\caption{Forbidden detour move for a virtual crossings of type $\alpha \not= 0$} \label{forbdet}
\end{figure}

Multi-virtual knot theory is a generalization of virtual knot theory in which a diagram may contain several distinct types of virtual crossings instead of a single kind. A multi-virtual knot diagram is an immersion of a circle in the plane whose double points are either classical crossings or virtual crossings labeled by a symbol $\alpha$, where $\alpha$ is simply a label taken from a fixed finite set (one may take it as an index set, e.g., integers, if desired). Equivalence is generated by the usual classical Reidemeister moves together with virtual Reidemeister–type moves that preserve the labels of virtual crossings. In particular, detour moves are allowed only when they respect the types (labels) of the virtual crossings involved. For more clarity see Figure~\ref{det} and \ref{forbdet}. This framework reduces to ordinary virtual knot theory when all virtual crossings share one type. The motivation for this generalization is that multiple virtual types encode additional combinatorial and topological data, which allows finer distinctions between knots and supports the construction of stronger invariants.
Kauffman developed this theory in the course of addressing the Penrose evaluation problem for colorings of arbitrary trivalent graphs.


\section{Multi-virtual braid group} 
\label{pure}

In this section, we recall (see \cite{Ka1}) definitions of the multi-virtual  and the symmetric multi-virtual braid groups.
The introduction of these groups is evidently motivated by multi-virtual
knot theory recalled above.

Let $k \geq 1$, 
$k$-{\it multiple virtual braid group} on $n$ strands, $M_kVB_n$ is a group which is generated by elements
$$
\sigma_{i}, \rho^{(\alpha)}_i, ~~i = 1, 2, \ldots, n-1,~\alpha = 0, 1, \ldots, k-1,
$$
and is defined by the following types of relations:

I. {\it Involutivity of generators}
$$
 \left( \rho^{(\alpha)}_i \right)^2 = e, ~~i = 1, 2, \ldots, n-1,~\alpha = 0, 1, \ldots, k-1.
$$

II. {\it For commutativity}

-- homogeneous
$$
\sigma_{i} \sigma_{j}= \sigma_j \sigma_{i}, \,\,  \rho_{i}^{(\alpha)} \rho_{j}^{(\alpha)} = \rho_{j}^{(\alpha)}   \rho_{i}^{(\alpha)}, ~~|i-j| \geq 2;
$$

-- mixed ($\alpha \not= \beta$)
$$
\sigma_{i} \rho_{j}^{(\alpha)}= \rho_{j}^{(\alpha)} \sigma_{i}, \,\,  \rho_{i}^{(\alpha)} \rho_{j}^{(\beta)} = \rho_{j}^{(\beta)}  \rho_{i}^{(\alpha)}, ~~|i-j| \geq 2.
$$

III. {\it Braid relations}

-- homogeneous
$$
\sigma_{i}  \sigma_{i+1}  \sigma_{i} = \sigma_{i+1}  \sigma_{i}  \sigma_{i+1}, \,\, \, \, 
\rho_{i}^{(\alpha)} \, \rho_{i+1}^{(\alpha)}  \, \rho_{i}^{(\alpha)} = \rho_{i+1}^{(\alpha)} \, \rho_{i}^{(\alpha)} \, \rho_{i+1}^{(\alpha)},~~
i=1, 2, \ldots, {n-2};
$$

-- mixed ($\beta = 1, 2, \ldots, k-1$)
$$
\sigma_{i}  \rho_{i+1}^{(0)}  \, \rho_{i}^{(0)} = \rho_{i+1}^{(0)} \, \rho_{i}^{(0)}\sigma_{i+1}, 
$$
$$
\rho_{i}^{(0)} \, \rho_{i+1}^{(0)}  \, \rho_{i}^{(\beta)} = \rho_{i+1}^{(\beta)} \, \rho_{i}^{(0)} \, \rho_{i+1}^{(0)},~~~~
i=1, 2, \ldots, {n-2}.
$$

\begin{remark}
Relations of the following forms are forbidden

$F1 \colon \sigma_{i}  \sigma_{i+1}  \, \rho_{i}^{(\alpha)} = \rho_{i+1}^{(\alpha)} \, \sigma_{i} \sigma_{i+1}$, $0 \leq \alpha \leq k-1$;

\medskip

$F2 \colon \sigma_{i+1}  \sigma_{i}  \, \rho_{i+1}^{(\alpha)} = \rho_{i}^{(\alpha)} \, \sigma_{i+1} \sigma_{i}$,  $0 \leq \alpha \leq k-1$;

\medskip

$F3 \colon \rho_{i}^{(0)} \, \rho_{i+1}^{(\beta)}  \, \rho_{i}^{(\beta)} = \rho_{i+1}^{(\beta)} \, \rho_{i}^{(\beta)} \, \rho_{i+1}^{(0)}$,~~
$ \rho_{i}^{(\gamma)} \, \rho_{i+1}^{(\beta)}  \, \rho_{i}^{(\beta)} = \rho_{i+1}^{(\beta)} \, \rho_{i}^{(\beta)} \, \rho_{i+1}^{(\gamma)}$,

\medskip

$\rho_{i}^{(\gamma)} \, \rho_{i+1}^{(\gamma)}  \, \rho_{i}^{(\beta)} = \rho_{i+1}^{(\beta)} \, \rho_{i}^{(\gamma)} \, \rho_{i+1}^{(\gamma)}$,~~
$0 < \gamma < \beta \leq k-1$.
\end{remark}

\begin{remark}
From the mixed braid relations  
it follows the relations
$$
\sigma_{i+1}  \rho_{i}^{(0)}  \, \rho_{i+1}^{(0)} = \rho_{i}^{(0)} \, \rho_{i+1}^{(0)}\sigma_{i} 
$$
and the  relations 
$$
\rho_{i+1}^{(0)} \, \rho_{i}^{(0)}  \, \rho_{i+1}^{(\beta)} = \rho_{i}^{(\beta)} \, \rho_{i+1}^{(0)} \, \rho_{i}^{(0)},~~~~
i=1, 2, \ldots, {n-1}.
$$
\end{remark}

If $k=1$, then $M_kVB_n$ is the virtual braid group which is denoted by $VB_n$. 

\medskip

\begin{proposition}
We have the following properties:
\begin{enumerate}
\item $VB_n \leq M_kVB_n$.
\item If $k \geq 2$, then $FVB_n \leq M_kVB_n$.
\item For $k\geq 1$, $M_kVB_n \leq M_{ k+1}VB_n$. 
\end{enumerate}
\end{proposition}

\begin{proof}
 \begin{itemize}
    \item[(1)] Define a map $\psi_1 \colon  M_kVB_n \to VB_n$ such that $\sigma_i \mapsto \sigma_i$, $\rho_i^{(0)} \mapsto \rho_i$ and $\rho_i^{(\alpha)} \mapsto e$. This map is a homomorphism.
Define a section of $\psi_1$, an inclusion map $\iota \colon  VB_n \hookrightarrow M_kVB_n$ such that $\psi_1 \circ \iota$ is the identity map on $VB_n$. The map $\iota \colon  VB_n \to M_kVB_n$ is a well defined homomorphism. Since the homomorphism $\psi_1$ has a section $\iota$ then the image $VB_n$ is a subgroup of $M_kVB_n$.
\item[(2)] Define a map $\psi_2 \colon   M_kVB_n \to FVB_n$ such that $\sigma_i \mapsto e$, $\rho_i^{(0)} \mapsto \rho_i$, $\rho_i^{(1)} \mapsto c_i$ and $\rho_i^{(\alpha)} \mapsto e$ for $\alpha\neq 0,1$. This map is a homomorphism.
Define a section of $\psi_2$, an inclusion map $\iota \colon  FVB_n \hookrightarrow M_kVB_n$ $c_i \mapsto \rho^{(1)}_i$ and $\rho_i \mapsto \rho_i^{(0)}$ such that $\psi_2 \circ \iota$ is the identity map on $FVB_n$. The map $\iota \colon  FVB_n \to M_kVB_n$ is a well defined homomorphism. Since the homomorphism $\psi_2$ has a section $\iota$ then the image $FVB_n$ is a subgroup of $M_kVB_n$.
\item[(3)] Define a map $\psi_3 \colon   M_{ k+1}VB_n \to M_kVB_n$ such that $\sigma_i \mapsto \sigma_i$, $\rho_i^{(\alpha)} \mapsto \rho^{(\alpha)}_i$ for $\alpha\neq k$, and $\rho_i^{(k)} \mapsto e$. This map is a homomorphism.
Define a section of $\psi_3$, an inclusion map $\iota: M_kVB_n \hookrightarrow M_{ k+1}VB_n$ such that $\psi_3 \circ \iota$ is the identity map on $M_kVB_n$. The map $\iota \colon  M_kVB_n \to M_{k+1}VB_n$ is a well defined homomorphism. Since the homomorphism~$\psi_3$ has a section $\iota$ then the image $M_kVB_n$ is a subgroup of $M_{ k+1}VB_n$.
\end{itemize}
\end{proof}

\subsection{Multi-virtual pure braid group}

Further we shall assume that $\rho_i = \rho_i^{(0)}$. 
Since, $S_n = \langle \rho_1, \rho_2, \ldots, \rho_{n-1} \rangle \leq VB_n$ and  $VB_n \leq M_kVB_n$, hence, $S_n$ is a subgroup of $M_kVB_n$.

Let us define a map 
$$
\varphi_{n,k} \colon M_kVB_n \to S_{n},~~\sigma_i \mapsto \rho_i,~~\rho_i^{(\alpha)} \mapsto \rho_i,~~i = 1, 2, \ldots, n-2,~\alpha = 0,1, 2, \ldots, k-1.  
$$ 
The kernel $\ker \varphi_{n,k}$ is the {\it  multi virtual pure braid group} $M_kVP_n$.

To find generators and  defining relations for  $M_kVP_n$, we use the Reidemeister-Schreier method (see, for example  \cite[Chapter 2.3]{MKS}). As a Schreier set of coset representation of $M_kVP_n$ in $M_kVB_n$ we take the same set $\Lambda_n$, which is used  in calculation of generators and relations for the group $VP_n$,
$$
\Lambda_n = \left\{ \prod\limits_{k=2}^n m_{k,j_k}~ |~ 1 \leq j_k
\leq k \right\}
$$
where $m_{kl}=\rho_{k-1}\rho_{k-2}\cdots \rho_l$ for $l<k$ and $m_{kl}=1$ in the other cases. 

Let us define elements that can be expressed  in terms of  the generators of $M_kVB_n$ by the formulas 
\begin{align*}
\lambda_{i,i+1}^{(0)} &= \rho_i \, \sigma_i^{-1},\\
\lambda_{i+1,i}^{(0)} &= \rho_i \, \lambda_{i,i+1}^{(0)} \, \rho_i = \sigma_i^{-1} \, \rho_i,
\end{align*}
\begin{align*}
\lambda_{i,i+1}^{(\beta)} &= \rho_i \, \rho_i^{(\beta)},~~1 \leq \beta \leq k-1,\\
\end{align*}
for  $i=1, 2, \ldots, n-1$, and
\begin{align*}
\lambda_{i,j}^{(0)} & =  \rho_{j-1} \, \rho_{j-2} \ldots \rho_{i+1} \, \lambda_{i,i+1}^{(0)} \, \rho_{i+1} \ldots \rho_{j-2} \, \rho_{j-1},\\
\lambda_{j,i}^{(0)} & =  \rho_{j-1} \, \rho_{j-2} \ldots \rho_{i+1} \, \lambda_{i+1,i}^{(0)} \, \rho_{i+1} \ldots \rho_{j-2} \, \rho_{j-1},\\
\lambda_{i,j}^{(\beta)} & = \rho_{j-1} \, \rho_{j-2} \ldots \rho_{i+1} \, \lambda_{i,i+1}^{(\beta)} \, \rho_{i+1} \ldots \rho_{j-2} \, \rho_{j-1}, 
\end{align*} 
for $1 \leq i < j-1 \leq n-1,~~1 \leq \beta \leq k-1$.
\begin{lemma}\label{gamma}
Let $a$ be an element of $\langle \rho_1, \rho_2, \ldots, \rho_{n-1} \rangle$ and $\bar{a}$ is its image in $S_n$ under the isomorphism $\rho_i \mapsto (i,i+1)$, $i = 1, 2, \ldots, n-1$, then for any generator $\lambda_{ij}^{(\alpha)}$ of $M_kVP_n$ the following holds
$$
a^{-1} \lambda_{ij}^{(\alpha)} a = \lambda_{(i)\bar{a}, (j)\bar{a}}^{(\alpha)},
$$
where $(k)\bar{a}$ is the image of $k$ under the action of the permutation $\bar{a}$.
\end{lemma}

\begin{theorem}\label{1purethm}
$M_kVP_n$  admits a presentation with the  generators 
$$
\lambda^{(0)}_{ij},~~1\leq i\neq j\leq n,
$$
$$
\lambda^{(\beta)}_{ij},~~1\leq i < j\leq n,~~1 \leq \beta \leq k-1,
$$
and the following relations:
\begin{align}
& \lambda_{ij}^{(\alpha)} \lambda_{kl}^{(\gamma)} =\lambda_{kl}^{(\gamma)} \lambda_{ij}^{(\alpha)} \label{1rel}, ~~0 \leq \alpha \leq \gamma \leq k-1,\\
&
\lambda_{ik}^{(0)} \lambda_{jk}^{(0)} \lambda_{ij}^{(0)} =\lambda_{ij}^{(0)} \lambda_{jk}^{(0)} \lambda_{ik}^{(0)},~~1 \leq i, j,  k \leq n,
\label{relation}\\
&
\lambda_{ik}^{(\beta)} \lambda_{jk}^{(\beta)} \lambda_{ij}^{(\beta)} =\lambda_{ij}^{(\beta)} \lambda_{jk}^{(\beta)} \lambda_{ik}^{(\beta)},~~1 \leq i < j < k \leq n, ~~1 \leq \beta \leq k-1
\label{relation1},
\end{align}
where distinct letters stand for distinct indices.

\end{theorem}

 \begin{proof}\label{main}
    Define the map $\Bar{} \colon M_kVB_n \to \Lambda_n$ which takes an element $w \in M_kVB_n$ to its representative $\overline{w}$ from $\Lambda_n$. In this case the element $w\overline{w}^{-1}$ belongs to $M_kVP_n$. By Theorem 2.7 of~\cite{MKS} the group $M_kVP_n$ is generated by
    $$
    s_{\lambda,a}=\lambda a \cdot (\overline{\lambda a})^{-1},~~~ \lambda \in \Lambda_n, ~~a \in \{\rho^{(\alpha)}_1, \ldots, \rho^{(\alpha)}_{n-1}, \sigma_1, \ldots, \sigma_{n-1}\},
    $$where $\alpha=0, 1, \ldots,k-1$. 
    
If $\alpha = 0$, then all   $s_{\lambda, \rho^{(0)}_i}=e$ and $s_{\lambda, \sigma_i}=\lambda (s_{e,\sigma_i} )\lambda^{-1}=\lambda (\sigma_i\rho_i)\lambda^{-1}=\lambda (\lambda_{i,i+1}^{-1})\lambda^{-1}$, which is equal to some $\lambda^{(0)}_{kl}$, by Lemma~\ref{form}.
   These calculations are done in Theorem 1~\cite{B}.
   
    Now, consider the generators such that $\alpha \neq 0$.
$$s_{\lambda,\rho^{(\alpha)}_i}=\lambda(s_{e,\rho^{(\alpha)}_i})\lambda^{-1}.$$
    Since $s_{e,\rho^{(\alpha)}_i}=\rho^{(\alpha)}_i\rho^{(0)}_i$, $s_{\lambda,\rho^{(\alpha)}_i}=\lambda(\rho^{(\alpha)}_i\rho^{(0)}_i)\lambda^{-1}$, which is equal to some $\lambda^{(\alpha)}_{kl}$ by Lemma~\ref{gamma}.
   Therefore, generators of the group $M_kVP_n$ are $$\lambda^{(0)}_{kl}, \text{ for } 1 \leq k\neq l \leq n \text{ and }
   \lambda^{(\beta)}_{ij},~~1\leq i < j\leq n,~~1 \leq \beta \leq k-1.$$

   To find the defining relations of $M_kVP_n$, we define a rewriting process $\tau$. It helps to rewrite a word $u$ to $\tau(u)$, where $u$ is written in the  generators of $M_kVB_n$ but represents an element of $M_kVP_n$ and $\tau(u)$ is a word written in  the generators of $M_kVP_n$. Let us associate to reduce word 
$$
u=a_1^{\epsilon_1}a_2^{\epsilon_2}\cdots a_v^{\epsilon_v}, ~~\epsilon_l=\pm 1, ~~a_l \in \{\sigma_1,\sigma_2, \ldots, \sigma_{n-1}, \rho^{(0)}_1,\rho^{(0)}_2, \ldots, \rho^{(0)}_{n-1},\ldots, \rho^{(k-1)}_1,\rho^{(k-1)}_2, \ldots, \rho^{(k-1)}_{n-1}\},
$$
 the word
   $$\tau(u)=s_{k_1,a_1}^{\epsilon_1}s_{k_2,a_2}^{\epsilon_2}\cdots s_{k_v,a_v}^{\epsilon_v},$$
   in the generators of $M_kVP_n$, where $k_j$ is the $(j-1)^{th}$ initial segment of the word $u$ if $\epsilon_j=1$, and a representative of the $j$-th initial segment of $u$ if $\epsilon_j=-1$.

   By Theorem 2.9 in~\cite{MKS}, the group $M_kVP_n$ is defined by the relations
   $$r_{\mu, \lambda}=\tau(\lambda r_\mu \lambda^{-1})=\lambda\tau( r_\mu) \lambda^{-1}, ~~~\lambda \in \Lambda_n,$$
   where $r_\mu$ is a defining relation of $M_kVB_n$.

Let us consider relations of $M_kVB_n$ stated in Section~\ref{pure}.
   Denote by $r_1=\left( \rho^{(\alpha)}_i \right)^2 = 1$. Then
   \begin{align*}
       r_{1,e} = \tau(r_1)& = s_{e,\rho^{(\alpha)}_i}s_{\rho_i, \rho^{(\alpha)}_i}\\
                         & = (\rho^{(\alpha)}_i\rho^{(0)}_i)(\rho^{(0)}_i\rho^{(\alpha)}_i\rho^{(0)}_i\rho^{(0)}_i)\\
                            & =e.
   \end{align*}
      Denote by $r_2=\sigma_{i} \sigma_{j} \sigma_{i}^{-1}\sigma_j^{-1}$, $|i-j| \geq 2$. Then
   \begin{align*}
       r_{2,e} = \tau(r_2)& = s_{e,\sigma_{i}}s_{\rho_i, \sigma_{j}}s^{-1}_{\rho_i\rho_j\rho_i,\sigma_{i}}s^{-1}_{\rho_i\rho_j\rho_i\rho_j, \sigma_{j}}\\
                         & = (\lambda^{(0)}_{i,i+1})^{-1}(\lambda^{(0)}_{j,j+1})^{-1}(\lambda^{(0)}_{i,i+1})(\lambda^{(0)}_{j,j+1}).
   \end{align*}
   The remaining $r_{2,\lambda}$, for all $\lambda \in \Lambda_n$, can be obtained from this relation using conjugation by $\lambda^{-1}$ and it gives the same relation, by Lemma~\ref{form}. We have obtained  $\lambda_{ij}^{(0)} \lambda_{kl}^{(0)} =\lambda_{kl}^{(0)} \lambda_{ij}^{(0)}$, where indices are different.

   Now, consider the next relation $r_3=\rho_{i}^{(\alpha)} \rho_{j}^{(\alpha)}\rho_{i}^{(\alpha)}\rho_{j}^{(\alpha)}$, $|i-j| \geq 2$.
   We have,  \begin{align*}
              r_{3,e} = \tau(r_3)& = s_{e,\rho_{i}^{(\alpha)}}s_{\rho_{i},\rho_{j}^{(\alpha)}}s_{\rho_{i}\rho_j,\rho_{i}^{(\alpha)}}s_{\rho_{i}\rho_j\rho_i,\rho_{j}^{(\alpha)}}\\
                         & = (\lambda_{i,i+1}^{(\alpha)})^{-1}(\lambda_{j,j+1}^{(\alpha)})^{-1}(\lambda_{i,i+1}^{(\alpha)})^{-1}(\lambda_{j,j+1}^{(\alpha)})^{-1}.
            \end{align*}
The remaining $r_{3,\lambda}$, for all  $\lambda \in \Lambda_n$. 
\begin{align*}
              r_{3,\lambda} = \tau(r_3)& = s_{\lambda,\rho_{i}^{(\alpha)}}s_{\lambda\rho_{i},\rho_{j}^{(\alpha)}}s_{\lambda\rho_{i}\rho_j,\rho_{i}^{(\alpha)}}s_{\lambda\rho_{i}\rho_j\rho_i,\rho_{j}^{(\alpha)}}\\
                         & = (\lambda(\lambda_{i,i+1}^{(\alpha)})^{-1}\lambda^{-1})(\lambda((\lambda_{j,j+1}^{(\alpha)})^{-1}\lambda^{-1})(\lambda((\lambda_{i,i+1}^{(\alpha)})^{-1}\lambda^{-1})(\lambda((\lambda_{j,j+1}^{(\alpha)})^{-1}\lambda^{-1})\\
                         &= \lambda(\lambda_{i,i+1}^{(\alpha)})^{-1}(\lambda_{j,j+1}^{(\alpha)})^{-1}(\lambda_{i,i+1}^{(\alpha)})^{-1}(\lambda_{j,j+1}^{(\alpha)})^{-1}\lambda^{-1}.
            \end{align*}

By Lemma~\ref{gamma}, we have obtained $\lambda_{ij}^{(\alpha)} \lambda_{kl}^{(\alpha)} =\lambda_{kl}^{(\alpha)} \lambda_{ij}^{(\alpha)}$, where indices are different.
The other commutative relations give rises to the relations: $\lambda_{ij}^{(\alpha)} \lambda_{kl}^{(\beta)} =\lambda_{kl}^{(\beta)} \lambda_{ij}^{(\alpha)}$, $0 \leq \alpha \leq \beta \leq k-1$. Combining all these relations, we obtain the relation~(\ref{1rel}).

Consider the next relation $r_4=\sigma_{i}  \sigma_{i+1}  \sigma_{i}\sigma^{-1}_{i+1}  \sigma^{-1}_{i} \sigma^{-1}_{i+1}$.
\begin{align*}
           r_{4,e} = \tau(r_4)& = s_{e,\sigma_i}s_{\rho_i,\sigma_{i+1}}s_{\rho_i\rho_{i+1},\sigma_i}s^{-1}_{\rho_{i+1}\rho_{i},\sigma_{i+1}}s_{\rho_{i+1}, \sigma_{i}} s^{-1}_{e,\sigma_{i+1}}\\
                         & = (\lambda^{(0)}_{i,i+1})^{-1}(\lambda^{(0)}_{i,i+2})^{-1}(\lambda^{(0)}_{i+1,i+2})^{-1}\lambda^{(0)}_{i,i+1}\lambda^{(0)}_{i,i+2}\lambda^{(0)}_{i+1,i+2}.
            \end{align*}
            The remaining relations $r_{4,\lambda}$, for all  $\lambda \in \Lambda_n$, can be obtained from previous relation using conjugation by $\lambda^{-1}$. By Lemma~\ref{form}, we obtain (\ref{relation}).

Let us consider the relation $r_5=\rho_{i}^{(\alpha)} \, \rho_{i+1}^{(\alpha)}  \, \rho_{i}^{(\alpha)} = \rho_{i+1}^{(\alpha)} \, \rho_{i}^{(\alpha)} \, \rho_{i+1}^{(\alpha)}$.

For $\alpha=0$, we get trivial relation. Now, solve for $\alpha\neq 0$.
\begin{align*}
              r_{5,e} = \tau(r_5)= & s_{e,\rho^{(\alpha)}_i}s_{\rho_i,\rho^{(\alpha)}_{i+1}}s_{\rho_i\rho_{i+1},\rho^{(\alpha)}_i}s^{-1}_{\rho_{i+1}\rho_{i},\rho^{(\alpha)}_{i+1}}s_{\rho_{i+1}, \rho^{(\alpha)}_{i}} s^{-1}_{e,\rho^{(\alpha)}_{i+1}}\\
                         & = (\lambda^{(\alpha)}_{i,i+1})^{-1}(\lambda^{(\alpha)}_{i,i+2})^{-1}(\lambda^{(\alpha)}_{i+1,i+2})^{-1}\lambda^{(\alpha)}_{i,i+1}\lambda^{(\alpha)}_{i,i+2}\lambda^{(\alpha)}_{i+1,i+2}.
            \end{align*}
Conjugating this relation by all representatives from $\Lambda_n$, by Lemma~\ref{gamma}, we obtain (\ref{relation1}).
Rest of the mixed relations of $M_kVB_n$ become trivial relations.
Therefore, the group $M_kVP_n$ is defined by the relations (\ref{1rel})-(\ref{relation1}).
\end{proof}

We have decomposition $M_kVB_n = M_kVP_n \rtimes S_n$.

\begin{corollary}\label{1purecoro}
We have the following properties:
\begin{enumerate}
\item $\langle \lambda_{ij}^{(0)} ~|~1 \leq i \not= j \leq n\rangle \cong VP_n$.
\item For any $\alpha \in \{ 1, 2, \ldots, k-1 \}$ holds $\langle \lambda_{ij}^{(\alpha)} ~|~1 \leq i \not= j \leq n\rangle \cong FVP_n$.
\item If $n=3$, then $M_kVP_3 \cong VP_3 * \underbrace{(FVP_3)*\cdots*(FVP_3)}_{(k-1)}$.
\end{enumerate}
\end{corollary}

{\begin{proof}
 \begin{itemize}
   
    \item[(1)] Define a endomorphism from $\phi_1 \colon \langle \lambda_{ij}^{(0)} ~|~1 \leq i \not= j \leq n\rangle \to VP_n$ such that $$\lambda_{ij}^{(0)} \to \lambda_{ij}.$$
    Define a section of $\phi_1$, an inclusion map $\iota \colon VP_n \hookrightarrow \langle \lambda_{ij}^{(0)} ~|~1 \leq i \not= j \leq n\rangle$ such that $\phi_1 \circ \iota$ is the identity map on $VP_n$. Since the homomorphism $\phi_1$ is bijective,  $\langle \lambda_{ij}^{(0)} ~|~1 \leq i \not= j \leq n\rangle \cong VP_n$.
    \item[(2)]  Let $\alpha \in \{ 1, 2, \ldots, k-1 \}$. Using the endomorphism $\phi_2 \colon \langle \lambda_{ij}^{(\alpha)} ~|~1 \leq i \not= j \leq n\rangle \to FVP_n$ such that $\lambda_{ij}^{(\alpha)} \to \lambda_{ij}$, as in the previous case, we get the need assertion.
    \item[(3)] It is easy to observe that $$M_kVP_3 = \langle \lambda_{ij}^{(0)} ~|~1 \leq i \not= j \leq 3\rangle *   \underbrace{\langle \lambda_{ij}^{(\alpha)} ~|~1 \leq i \not= j \leq 3\rangle*\cdots*\langle \lambda_{ij}^{(\alpha)} ~|~1 \leq i \not= j \leq 3\rangle}_{(k-1)}.$$ Using (2) and (3), we obtain the result $M_kVP_3 \cong VP_3  * \underbrace{(FVP_3)*\cdots*(FVP_3)}_{(k-1)}$.
\end{itemize}
\end{proof}

\begin{example}
To illustrate definitions and results of this  sections, we consider 2-virtual 3-strand braid group. As a rule, 3-strand braid groups are different from the braid groups with larger number of strands and are therefore studied separately (see \cite{BK}, \cite{BMVW}, \cite{GKM}, \cite{K-4}).

By the definition, the group $MVB_3 = M_2VB_3$ and it is generated by 6 elements:
$$
\sigma_1,~~ \sigma_2,~~\rho_1,~~\rho_2,~~\tau_1,~~\tau_2,
$$
and is defined by relations
$$
\rho_1^2 = \rho_2^2 = \tau_1^2 = \tau_2^{2} = e,
$$
$$
\sigma_{1}  \sigma_{2}  \sigma_{1} = \sigma_{2}  \sigma_{1}  \sigma_{2},~~\rho_1 \rho_2 \rho_1 = \rho_2 \rho_1 \rho_2,~~\tau_1 \tau_2 \tau_1 = \tau_2 \tau_1 \tau_2,
$$
$$
\sigma_{1}  \rho_{2}  \rho_{1} =  \rho_{2} \rho_{1} \sigma_{2},~~\tau_1 \rho_2 \rho_1 = \rho_2 \rho_1 \tau_2.
$$
It is not difficult to see that  $\langle \sigma_1,~~ \sigma_2,~~\rho_1,~~\rho_2 \rangle \cong VB_3$ is the virtual braid group and  $\langle \sigma_1,~~ \sigma_2,~~\tau_1,~~\tau_2 \rangle \cong FVB_3$ is the flat virtual braid group.

The kernel $\ker \varphi_{3}$ of the homomorphism $\varphi_{3} \colon MVB_3 \to S_{3}$,
 which is defined   on the generators:
$$
\sigma_i \mapsto \rho_i,~~\tau_i \mapsto \rho_i,~~\rho_i \mapsto \rho_i,~~i = 1, 2.  
$$ 
The $\ker \varphi_{3}$ is {\it  2-virtual  pure  braid group} $MVP_3$.
This group is generated by elements
$$
\lambda_{12} = \rho_1 \sigma_1^{-1},~~\lambda_{21} =\sigma_1^{-1}  \rho_1,~~\lambda_{23} = \rho_2 \sigma_2^{-1},~~\lambda_{32} =\sigma_2^{-1}  \rho_2,
$$
$$
\lambda_{13} = \rho_2 (\rho_1 \sigma_1^{-1}) \rho_2 = \rho_2 \lambda_{12} \rho_2,~~\lambda_{31} = \rho_2 (\sigma_1^{-1}  \rho_1) \rho_2 = \rho_2 \lambda_{21} \rho_2.
$$
$$
\mu_{12} = \rho_1 \tau_1,~~\mu_{23} = \rho_2 \tau_2,~~
\mu_{13} = \rho_2 (\rho_1 \tau_1) \rho_2 = \rho_2 \mu_{12} \rho_2,
$$

and is defined by the relations
$$
\lambda_{12} \lambda_{32} \lambda_{13} = \lambda_{13} \lambda_{32} \lambda_{12},
~~~
\lambda_{21} \lambda_{31} \lambda_{23} = \lambda_{23} \lambda_{31} \lambda_{21},
$$
$$
\lambda_{13} \lambda_{23} \lambda_{12} = \lambda_{12} \lambda_{23} \lambda_{13},
~~~
\lambda_{31} \lambda_{21} \lambda_{32} = \lambda_{32} \lambda_{21} \lambda_{31},
$$
$$
\lambda_{23} \lambda_{13} \lambda_{21} = \lambda_{21} \lambda_{13} \lambda_{23},
~~~
\lambda_{32} \lambda_{12} \lambda_{31} = \lambda_{31} \lambda_{12} \lambda_{32},
$$
$$
\mu_{12} \mu_{13} \mu_{23} = \mu_{23} \mu_{13} \mu_{12}.
$$

\end{example}

\subsection{Multi virtual semi--pure braid group}

Define a  map
$$
\psi_{n,k} \colon M_kVB_n \to S_{n},~~\sigma_i \mapsto e,~~\rho_i^{(\alpha)} \mapsto \rho_i,~~i = 1, 2, \ldots, n-2,~\alpha = 0, 1, 2, \ldots, k-1.  
$$ 
The kernel $\ker \psi_{n,k}$ is said to be a {\it multi  semi-pure virtual braid group} $M_kVH_n$.
\begin{lemma}\label{hl}
Let $a$ be an element of $\langle \rho_1, \rho_2, \ldots, \rho_{n-1} \rangle$ and $\bar{a}$ is its image in $S_n$ under the isomorphism $\rho_i \mapsto (i,i+1)$, $i = 1, 2, \ldots, n-1$, then for any generator $x_{ij}^{(\alpha)}$ of $M_kVH_n$ the following holds
$$
a^{-1} x_{ij}^{(\alpha)} a = x_{(i)\bar{a}, (j)\bar{a}}^{(\alpha)},
$$
where $(k)\bar{a}$ is the image of $k$ under the action of the permutation $\bar{a}$ and we use the relations
$$
x_{ji}^{(\beta)} = \left(x_{ij}^{(\beta)}\right)^{-1},~~~~1 \leq \beta \leq k-1.
$$
\end{lemma}
\begin{theorem}
$M_kVH_n$  admits a presentation with the  generators 
$$
x^{(0)}_{ij},~~1\leq i\neq j\leq n,
$$
$$
x^{(\beta)}_{ij},~~1\leq i < j\leq n,~~1 \leq \beta \leq k-1,
$$
and the following relations:
\begin{align}
& x_{ij}^{(\alpha)} x_{kl}^{(\gamma)} =x_{kl}^{(\gamma)} x_{ij}^{(\alpha)} \label{rel}, ~~0 \leq \alpha \leq \gamma \leq k-1,\\
&
x_{ik}^{(0)} x_{kj}^{(0)} x_{ik}^{(0)} =x_{kj}^{(0)} x_{ik}^{(0)} x_{kj}^{(0)},~~1 \leq i, j,  k \leq n,
\label{relation}\\
&
x_{ik}^{(\beta)} x_{jk}^{(\beta)} x_{ij}^{(\beta)} = x_{ij}^{(\beta)} x_{jk}^{(\beta)} x_{ik}^{(\beta)},~~1 \leq i < j < k \leq n, ~~1 \leq \beta \leq k-1
\label{relation1},
\end{align}
where distinct letters stand for distinct indices.
The generators of $M_kVH_n$ can be expressed  in terms of  the generators of $M_kVB_n$ by the formulas 
\begin{align*}
x_{i,i+1}^{(0)} &= \sigma_i^{-1},\\
x_{i+1,i}^{(0)} &= \rho_i \, x_{i,i+1}^{(0)} \, \rho_i =  \rho_i \sigma_i^{-1} \, \rho_i,
\end{align*}
\begin{align*}
x_{i,i+1}^{(\beta)} &= \rho_i \, \rho_i^{(\beta)},~~1 \leq \beta \leq k-1,\\
\end{align*}
for  $i=1, 2, \ldots, n-1$, and
\begin{align*}
x_{i,j}^{(\alpha)} & = \rho_{j-1} \, \rho_{j-2} \ldots \rho_{i+1} \, x_{i,i+1}^{(\alpha)} \, \rho_{i+1} \ldots \rho_{j-2} \, \rho_{j-1}, \\ 
x_{j,i}^{(0)} & =  \rho_{j-1} \, \rho_{j-2} \ldots \rho_{i+1} \, x_{i+1,i}^{(0)} \, \rho_{i+1} \ldots \rho_{j-2} \, \rho_{j-1},
\end{align*} 
for $1 \leq i < j-1 \leq n-1,~~0 \leq \alpha \leq k-1$.
\end{theorem}

\begin{proof}
 It is similar to the proof of Theorem~\ref{1purethm}.
Define the map $\Bar{} \colon M_kVB_n \to \Lambda_n$ which takes an element $w \in M_kVB_n$ to its representative $\overline{w}$ from $\Lambda_n$. In this case the element $w\overline{w}^{-1}$ belongs to $M_kVH_n$. By Theorem 2.7 of~\cite{MKS} the group $M_kVH_n$ is generated by
    $$
    s_{\lambda,a}=\lambda a \cdot (\overline{\lambda a})^{-1},~~~ \lambda \in \Lambda_n, ~~a \in \{\rho^{(\alpha)}_1, \ldots, \rho^{(\alpha)}_{n-1}, \sigma_1, \ldots, \sigma_{n-1}\},
    $$where $\alpha=0, 1, \ldots,k-1$. 
    
If $\alpha = 0$, then all   $s_{\lambda, \rho^{(0)}_i}=e$ and $s_{\lambda, \sigma_i}=\lambda (s_{e,\sigma_i} )\lambda^{-1}=\lambda (\sigma_i)\lambda^{-1}=\lambda (x^{(0)}_{i,i+1})\lambda^{-1}$, which is equal to some $x^{(0)}_{kl}$, by Lemma~\ref{form1}.
   These calculations are done in Proposition 17~\cite{BB}.
   
    Now, consider the generators such that $\alpha \neq 0$.
$$s_{\lambda,\rho^{(\alpha)}_i}=\lambda(s_{e,\rho^{(\alpha)}_i})\lambda^{-1}.$$
    Since $s_{e,\rho^{(\alpha)}_i}=\rho^{(\alpha)}_i\rho^{(0)}_i$, $s_{\lambda,\rho^{(\alpha)}_i}=\lambda(\rho^{(\alpha)}_i\rho^{(0)}_i)\lambda^{-1}$, which is equal to some $x^{(\alpha)}_{kl}$ by Lemma~\ref{hl}.
   Therefore, generators of the group $M_kVH_n$ are $$x^{(0)}_{kl}, \text{ for } 1 \leq k\neq l \leq n \text{ and }
   x^{(\beta)}_{ij},~~1\leq i < j\leq n,~~1 \leq \beta \leq k-1.$$

   To find the defining relations of $M_kVH_n$, we define a rewriting process $\tau$. It helps to rewrite a word $u$ to $\tau(u)$, where $u$ is written in the  generators of $M_kVB_n$ but represents an element of $M_kVH_n$ and $\tau(u)$ is a word written in  the generators of $M_kVH_n$. Let us associate to reduce word 
$$
u=a_1^{\epsilon_1}a_2^{\epsilon_2}\cdots a_v^{\epsilon_v}, ~~\epsilon_l=\pm 1, ~~a_l \in \{\sigma_1,\sigma_2, \ldots, \sigma_{n-1}, \rho^{(0)}_1,\rho^{(0)}_2, \ldots, \rho^{(0)}_{n-1},\ldots, \rho^{(k-1)}_1,\rho^{(k-1)}_2, \ldots, \rho^{(k-1)}_{n-1}\},
$$
 the word
   $$\tau(u)=s_{k_1,a_1}^{\epsilon_1}s_{k_2,a_2}^{\epsilon_2}\cdots s_{k_v,a_v}^{\epsilon_v},$$
   in the generators of $M_kVH_n$, where $k_j$ is the $(j-1)^{th}$ initial segment of the word $u$ if $\epsilon_j=1$, and a representative of the $j$-th initial segment of $u$ if $\epsilon_j=-1$.

   By Theorem 2.9 in~\cite{MKS}, the group $M_kVH_n$ is defined by the relations
   $$r_{\mu, \lambda}=\tau(\lambda r_\mu \lambda^{-1})=\lambda\tau( r_\mu) \lambda^{-1}, ~~~\lambda \in \Lambda_n,$$
   where $r_\mu$ is a defining relation of $M_kVB_n$.

Let us consider relations of $M_kVB_n$ stated in Section~\ref{pure}.

Consider the $r=\sigma_{i} \sigma_{j} \sigma_{i}^{-1}\sigma_j^{-1}$, $|i-j| \geq 2$. Then
   \begin{align*}
       \tau(r)& = s_{e,\sigma_{i}}s_{\rho_i, \sigma_{j}}s^{-1}_{\rho_i\rho_j\rho_i,\sigma_{i}}s^{-1}_{\rho_i\rho_j\rho_i\rho_j, \sigma_{j}}\\
                         & = (x^{(0)}_{i,i+1})^{-1}(x^{(0)}_{j,j+1})^{-1}(x^{(0)}_{i,i+1})(x^{(0)}_{j,j+1}).
   \end{align*}

The remaining relations $\lambda \tau(r)\lambda^{-1}$, for all $\lambda \in \Lambda_n$, can be obtained from previous relation using conjugation by $\lambda^{-1}$ and it gives the same relation, by Lemma~\ref{form1}. We have obtained  $x_{ij}^{(0)} x_{kl}^{(0)} =x_{kl}^{(0)} x_{ij}^{(0)}$, where indices are different.

We can work on other relations of $M_kVB_n$ similarly as done in Theorem~\ref{1purethm}.
\end{proof}

\begin{corollary} 
We have the following properties:
\begin{enumerate}
\item $M_kVB_n = M_kVH_n \rtimes S_n$.
\item $\langle x_{ij}^{(0)} ~|~1 \leq i \not= j \leq n\rangle \cong VH_n$.
\item For any $\beta \in \{ 1, 2, \ldots, k-1 \}$ holds $\langle x_{ij}^{(\beta)} ~|~1 \leq i \not= j \leq n\rangle \cong FVP_n$.
\item If $n=3$, then $M_kVH_3 \cong VH_3  * \underbrace{(FVP_3)*\cdots*(FVP_3)}_{(k-1)}$.
\end{enumerate}
\end{corollary}

{\begin{proof}
\begin{itemize}
    \item[(1)] This holds true based on the definition of $M_kVH_n$. 
    \item[(2)] Define a endomorphism from $\phi'_1 \colon \langle x_{ij}^{(0)} ~|~1 \leq i \not= j \leq n\rangle \to VH_n$ such that $$x_{ij}^{(0)} \to x_{ij}.$$
    Define a section of $\phi'_1$, an inclusion map $\iota \colon VH_n \hookrightarrow \langle x_{ij}^{(0)} ~|~1 \leq i \not= j \leq n\rangle$ such that $\phi'_1 \circ \iota$ is the identity map on $VH_n$. Since the homomorphism $\phi'_1$ is bijective. Hence, $\langle x_{ij}^{(0)} ~|~1 \leq i \not= j \leq n\rangle \cong VH_n$.
    \item[(3)] Similar to the above proof. Let $\beta \in \{ 1, 2, \ldots, k-1 \}$. Using the endomorphism $\phi'_2 \colon \langle x_{ij}^{(\beta)} ~|~1 \leq i \not= j \leq n\rangle \to FVP_n$ such that $x_{ij}^{(\beta)} \to x_{ij}$.
    \item[(4)] It is easy to observe that $$M_kVH_3 = \langle x_{ij}^{(0)} ~|~1 \leq i \not= j \leq 3\rangle *   \underbrace{\langle x_{ij}^{(\beta)} ~|~1 \leq i \not= j \leq 3\rangle*\cdots*\langle x_{ij}^{(\beta)} ~|~1 \leq i \not= j \leq 3\rangle}_{(k-1)}.$$ Using (2) and (3), we obtain the result $M_kVH_3 \cong VH_3  * \underbrace{(FVP_3)*\cdots*(FVP_3)}_{(k-1)}$.
\end{itemize}
\end{proof}

We have decomposition $M_kVB_n = M_kVH_n \rtimes S_n$ and $S_n$ acts on $M_kVH_n$ by the rules.

\newpage

\subsection{Symmetric multi-virtual braid group}

In the definition of $M_kVB_n$, virtual crossings need not detour over each other. The detour moves consist of only one type of virtual crossing.
Also, L. Kauffman \cite{Ka1}  suggested to study analogous of this group in which   are  allowed  detour moves  for any types of virtual crossings.

Let $k \geq 2$. 
{\it Symmetric 
 $k$-multiple virtual braid group} $\widetilde{M_kVB}_n$ is a group which is quotient of $M_kVB_n$ by the forbidden relations $F3$.

But in $\widetilde{M_kVB}_n$ the following relations are also forbidden:

\medskip

$F1 \colon \sigma_{i}  \sigma_{i+1}  \, \rho_{i}^{(\alpha)} = \rho_{i+1}^{(\alpha)} \, \sigma_{i} \sigma_{i+1}$, $0 \leq \alpha \leq k-1$;

\medskip

$F2 \colon \sigma_{i+1}  \sigma_{i}  \, \rho_{i+1}^{(\alpha)} = \rho_{i}^{(\alpha)} \, \sigma_{i+1} \sigma_{i}$,  $0 \leq \alpha \leq k-1$;

\vspace{0.5cm}
Hence, $\widetilde{M_kVB}_n$ is a group which is generated by elements
$$
\sigma_{i}, \rho^{(\alpha)}_i, ~~i = 1, 2, \ldots, n-1,~\alpha = 0, 1, \ldots, k-1,
$$
and is defined by the following types of relations:

\vspace{0.5cm}

I. {\it Involutivity of generators}
\medskip
$$
 \left( \rho^{(\alpha)}_i \right)^2 = 1, ~~i = 1, 2, \ldots, n-1,~\alpha = 0, 1, \ldots, k-1.
$$

\vspace{0.5cm}
II. {\it For commutativity}

\medskip

-- homogeneous
$$
\sigma_{i} \sigma_{j}= \sigma_j \sigma_{i}, \,\,  \rho_{i}^{(\alpha)} \rho_{j}^{(\alpha)} = \rho_{j}^{(\alpha)}   \rho_{i}^{(\alpha)}, ~~|i-j| \geq 2;
$$

\medskip

-- mixed ($\alpha \not= \beta$)
$$
\sigma_{i} \rho_{j}^{(\alpha)}= \rho_{j}^{(\alpha)} \sigma_{i}, \,\,  \rho_{i}^{(\alpha)} \rho_{j}^{(\beta)} = \rho_{j}^{(\beta)}  \rho_{i}^{(\alpha)}, ~~|i-j| \geq 2.
$$
\vspace{0.5cm}
III. {\it Braid relations}
\medskip

-- homogeneous
$$
\sigma_{i}  \sigma_{i+1}  \sigma_{i} = \sigma_{i+1}  \sigma_{i}  \sigma_{i+1}, \,\, \, \, 
\rho_{i}^{(\alpha)} \, \rho_{i+1}^{(\alpha)}  \, \rho_{i}^{(\alpha)} = \rho_{i+1}^{(\alpha)} \, \rho_{i}^{(\alpha)} \, \rho_{i+1}^{(\alpha)},~~
i=1, 2, \ldots, {n-2};
$$
\medskip

-- mixed ($\alpha < \beta$)
$$
\sigma_{i}  \rho_{i+1}^{(\alpha)}  \, \rho_{i}^{(\alpha)} = \rho_{i+1}^{(\alpha)} \, \rho_{i}^{(\alpha)}\sigma_{i+1}, 
$$
\vspace{0.3cm}
$$
\rho_{i}^{(\alpha)} \, \rho_{i+1}^{(\alpha)}  \, \rho_{i}^{(\beta)} = \rho_{i+1}^{(\beta)} \, \rho_{i}^{(\alpha)} \, \rho_{i+1}^{(\alpha)},~~~~
\rho_{i}^{(\alpha)} \, \rho_{i+1}^{(\beta)}  \, \rho_{i}^{(\beta)} = \rho_{i+1}^{(\beta)} \, \rho_{i}^{(\beta)} \, \rho_{i+1}^{(\alpha)},~~
i=1, 2, \ldots, {n-2}.
$$

\newpage

\begin{example}
For $M_2VB_n$, follow by Kauffman, we will use the following set of generators 
$$\left\{ \sigma_{i},  \rho_{i},\  \tau_{i}  |  1\leq i \leq n-1  \right\}:$$

\begin{figure}[ht]
\includegraphics[totalheight=4cm]{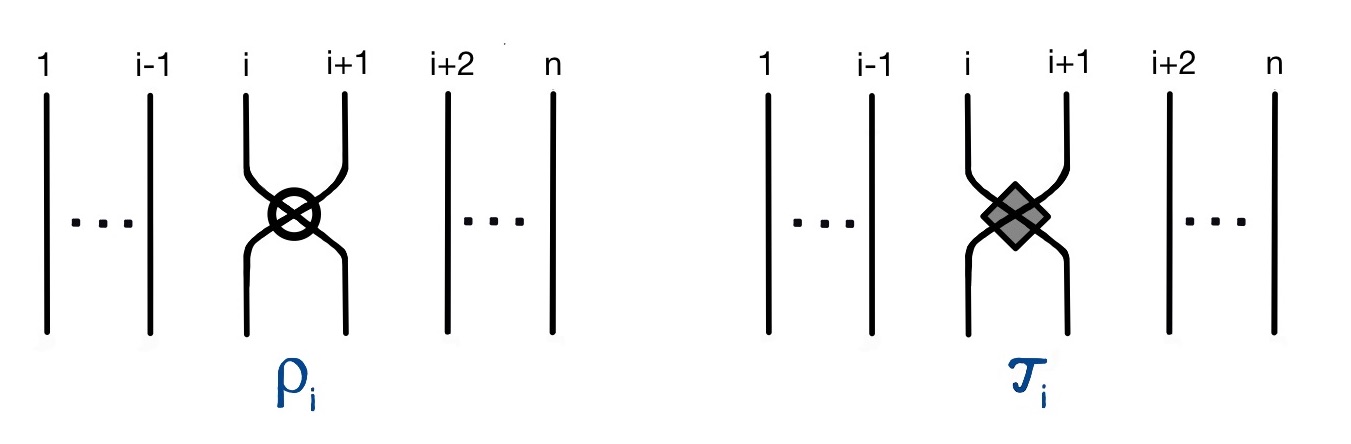}
\begin{flushleft}

and the defining relations:\\
\vspace{0.5cm}

$(1)  \,\,\ \sigma_{i} \sigma_{i}^{-1}= \sigma_i^{-1} \sigma_{i}= 1_n;$ 
\end{flushleft}
\includegraphics[totalheight=4cm]{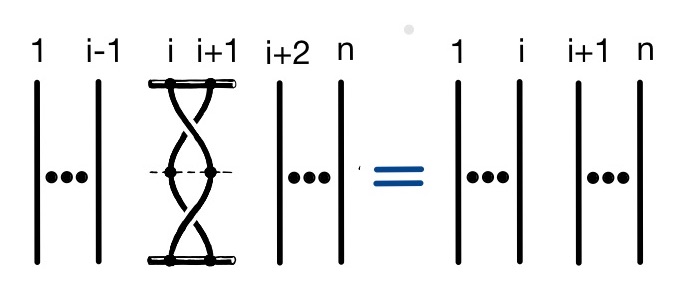}
\begin{flushleft}

\vspace{0.5cm}

$(2)  \,\,\ \sigma_{i} \sigma_{j}= \sigma_j \sigma_{i}, \,\,  |i-j| \geq 2;$
\end{flushleft}
\begin{flushleft}

\vspace{0.5cm}

$(3)  \,\,\ \sigma_{i}  \sigma_{i+1}  \sigma_{i} = \sigma_{i+1}  \sigma_{i}  \sigma_{i+1}, \,\,  i=1, 2, \ldots, {n-1};$ 
\end{flushleft}
\begin{flushleft}

\vspace{0.5cm}

$(4) \,\, \rho_{i}^{2}=1_{n},  \,\,  i=1, 2, \ldots, {n-1};$
\end{flushleft}

\includegraphics[totalheight=4cm]{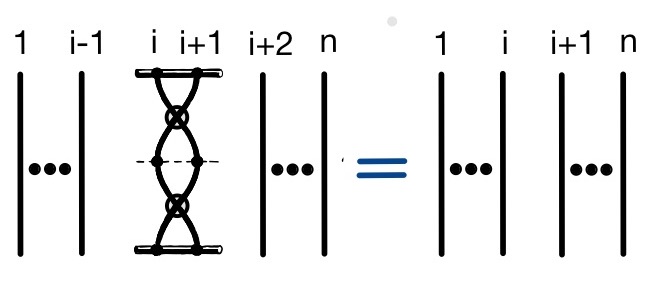}
\end{figure}

\begin{figure}
\begin{flushleft}

$(5) \,\, \tau_{i}^{2}=1_{n}, \,\,  i=1, 2, \ldots, {n-1};$
\end{flushleft}

\includegraphics[totalheight=3cm]{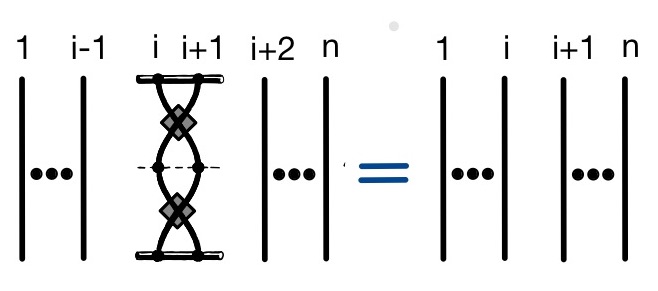}
\begin{flushleft}
\vspace{0.3cm}

(6) \,\, $\rho_i \tau_j=\tau_j \rho_i, \,\,  |i-j| \geq 2;$
\end{flushleft}
\includegraphics[totalheight=3.3cm]{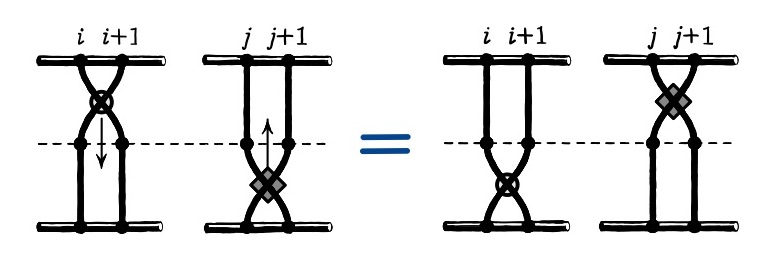}
\begin{flushleft}

\vspace{0.3cm}

$(7)  \,\,\ \rho_{i} \rho_{j}= \rho_j \rho_{i}, \,\,  |i-j| \geq 2;$
\end{flushleft}
\begin{flushleft}

\vspace{0.3cm}

$(8)  \,\,\ \tau_{i} \tau_{j}= \tau_j \tau_{i}, \,\,  |i-j| \geq 2;$
\end{flushleft}
\begin{flushleft}

\vspace{0.3cm}

$(9)  \,\,\ \rho_{i}  \rho_{i+1}  \rho_{i} = \rho_{i+1}  \rho_{i}  \rho_{i+1}, \,\,  i=1, 2, \ldots, {n-1};$ 
\end{flushleft}
\begin{flushleft}

\vspace{0.3cm}

$(10)  \,\,\ \tau_{i}  \tau_{i+1}  \tau_{i} = \tau_{i+1}  \tau_{i}  \tau_{i+1}, \,\,  i=1, 2, \ldots, {n-1};$ 
\end{flushleft}

\vspace{0.3cm}

\begin{flushleft}
$(11) \,\, \rho_{i} \tau_{i+1} \rho_{i} = \rho_{i+1}  \tau_{i} \rho_{i+1}, \,\,  i=1, 2, \ldots, {n-1};$
\end{flushleft}
\includegraphics[totalheight=4.2cm]{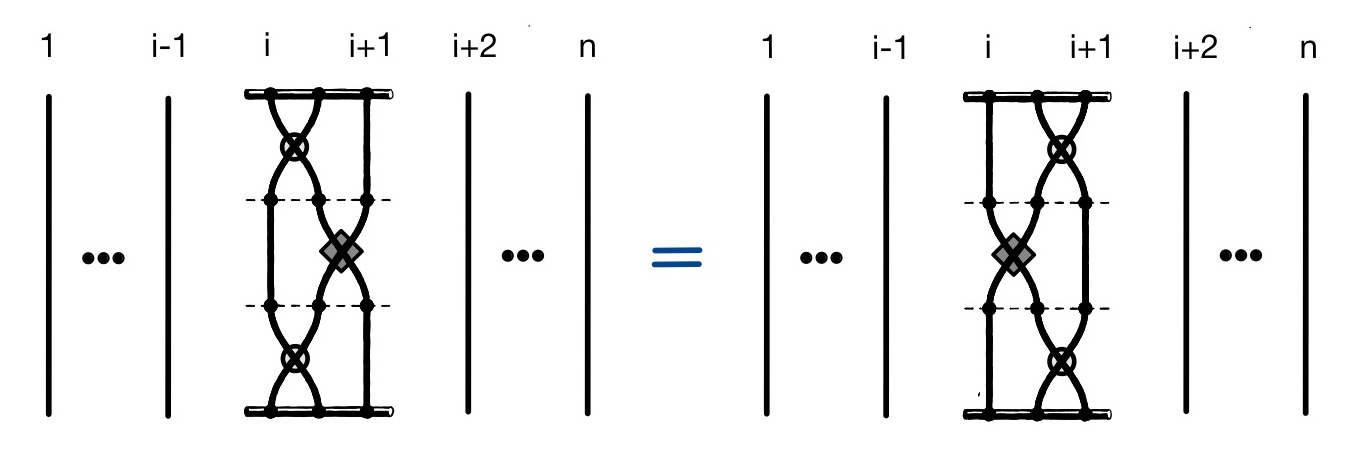}
\begin{flushleft}

\vspace{0.5cm}

$(12) \,\, \rho_{i} \sigma_{i+1} \rho_{i}=  \rho_{i+1} \sigma_{i} \rho_{i+1}, \,\,  i=1, 2, \ldots, {n-1};$
\end{flushleft}
\includegraphics[totalheight=4.2cm]{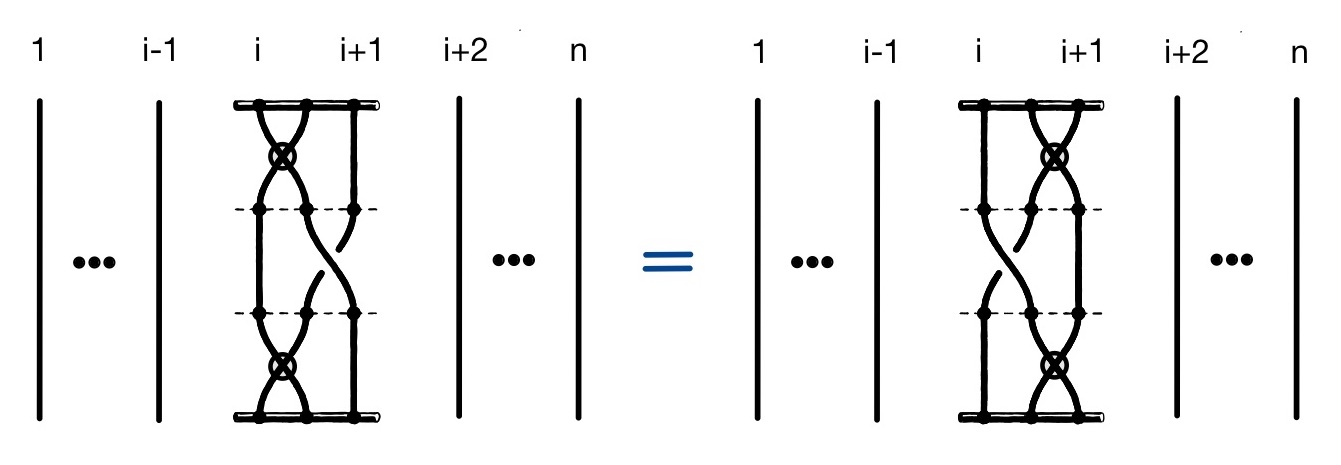}
\end{figure}
\begin{figure}
\begin{flushleft}
$(13) \,\, \rho_{i} \rho_{i+1} \sigma_{i}=  \sigma_{i+1} \rho_{i} \rho_{i+1} , \,\,  i=1, 2, \ldots, {n-1};$
\end{flushleft}
\end{figure}
\begin{figure}
\begin{flushleft}
$(14) \,\, \rho_{i} \ \rho_{i+1} \tau_{i}= \tau_{i+1}  \rho_{i} \rho_{i+1} , \,\,  i=1, 2, \ldots, {n-1}.$
\end{flushleft}
\end{figure}
\end{example}

\pagebreak

\subsection{Symmetric multi-virtual pure braid group}

Let us define a map
$$
\phi_{n,k} \colon ~~ \widetilde{M_kVB_n} \to S_{n},~~\sigma_i \mapsto \rho_i,~~\rho_i^{(\alpha)} \mapsto \rho_i,~~i = 1, 2, \ldots, n-2,~\alpha = 0, 1, 2, \ldots, k-1.  
$$ 
The kernel $\ker \phi_{n,k}$ is the {\it  symmetric multi-virtual pure  braid group} $\widetilde{M_kVP_n}$.

Using the fact that $\widetilde{M_kVB}_n$
 is the quotient of $M_kVB_n$ by the relations $$
\sigma_{i}  \rho_{i+1}^{(\alpha)}  \, \rho_{i}^{(\alpha)} = \rho_{i+1}^{(\alpha)} \, \rho_{i}^{(\alpha)}\sigma_{i+1}, ~~~~\text{ and }
~~~~\rho_{i}^{(\alpha)} \, \rho_{i+1}^{(\alpha)}  \, \rho_{i}^{(\beta)} = \rho_{i+1}^{(\beta)} \, \rho_{i}^{(\alpha)} \, \rho_{i+1}^{(\alpha)},$$ 
let us find a presentation for $\widetilde{M_kVP}_n$. 

From the relation $\sigma_{i}  \rho_{i+1}^{(\alpha)}  \, \rho_{i}^{(\alpha)} = \rho_{i+1}^{(\alpha)} \, \rho_{i}^{(\alpha)}\sigma_{i+1} \text{, for } \alpha\neq 0,$
follows the  relation:
$$\lambda_{ij}^{(0)}\lambda_{ik}^{(\alpha)}\lambda_{jk}^{(\alpha)}=\lambda_{jk}^{(\alpha)}\lambda_{ik}^{(\alpha)}\lambda_{ij}^{(0)}.$$

From the relation $\rho_{i}^{(\alpha)} \, \rho_{i+1}^{(\alpha)}  \, \rho_{i}^{(\beta)} = \rho_{i+1}^{(\beta)} \, \rho_{i}^{(\alpha)} \, \rho_{i+1}^{(\alpha)}$, for $0<\alpha\neq\beta$,
we obtain relation:
$$\lambda_{ij}^{(\alpha)}\lambda_{ik}^{(\alpha)}\lambda_{jk}^{(\beta)}=\lambda_{jk}^{(\beta)}\lambda_{ik}^{(\alpha)}\lambda_{ij}^{(\alpha)}.$$

From the relation $\rho_{i}^{(\alpha)} \, \rho_{i+1}^{(\alpha)}  \, \rho_{i}^{(0)} = \rho_{i+1}^{(0)} \, \rho_{i}^{(\alpha)} \, \rho_{i+1}^{(\alpha)}$, for $\alpha\neq 0$,
we obtain relation:
$$\lambda_{ij}^{(\alpha)}\lambda_{jk}^{(\alpha)}=\lambda_{jk}^{(\alpha)}\lambda_{ij}^{(\alpha)}.$$
This last relation is already in $M_kVP_n$.

By the above calculations, we obtain

\begin{proposition}
$\widetilde{M_kVP}_n$ is the quotient of $M_kVP_n$ by the relations
$$
\lambda_{ij}^{(0)}\lambda_{ik}^{(\alpha)}\lambda_{jk}^{(\alpha)}=\lambda_{jk}^{(\alpha)}\lambda_{ik}^{(\alpha)}\lambda_{ij}^{(0)},~~\alpha\neq 0,
$$
$$
\lambda_{ij}^{(\alpha)}\lambda_{ik}^{(\alpha)}\lambda_{jk}^{(\beta)}=\lambda_{jk}^{(\beta)}\lambda_{ik}^{(\alpha)}\lambda_{ij}^{(\alpha)},~~0<\alpha<\beta.
$$
\end{proposition}

\subsection{Symmetric multi-virtual semi-pure braid group}

Let us define a map 
$$
\psi_{n,k} \colon ~~ \widetilde{M_kVB_n} \to S_{n},~~\sigma_i \mapsto e,~~\rho_i^{(\alpha)} \mapsto \rho_i,~~i = 1, 2, \ldots, n-2,~\alpha = 0, 1, 2, \ldots, k-1.  
$$ 
The kernel $\ker \psi_{n,k}$ is the {\it  symmetric multi semi-pure virtual braid group} $\widetilde{M_kVH_n}$.

Let us find a presentation for $\widetilde{M_kVH}_n$. Since $\widetilde{M_kVB}_n$ it is the quotient $M_kVB_n$ by the relations $$
\sigma_{i}  \rho_{i+1}^{(\alpha)}  \, \rho_{i}^{(\alpha)} = \rho_{i+1}^{(\alpha)} \, \rho_{i}^{(\alpha)}\sigma_{i+1}, ~~~~\text{ and }
~~~~\rho_{i}^{(\alpha)} \, \rho_{i+1}^{(\alpha)}  \, \rho_{i}^{(\beta)} = \rho_{i+1}^{(\beta)} \, \rho_{i}^{(\alpha)} \, \rho_{i+1}^{(\alpha)},$$ 
from the relation $\sigma_{i}  \rho_{i+1}^{(\alpha)}  \, \rho_{i}^{(\alpha)} = \rho_{i+1}^{(\alpha)} \, \rho_{i}^{(\alpha)}\sigma_{i+1} \text{, for } \alpha\neq 0,$
we obtain the relation:
$$x_{ij}^{(0)}x_{ik}^{(\alpha)}x_{jk}^{(\alpha)}=x_{ik}^{(\alpha)}x_{jk}^{(\alpha)}x_{ij}^{(0)};$$
from the relation $\rho_{i}^{(\alpha)} \, \rho_{i+1}^{(\alpha)}  \, \rho_{i}^{(\beta)} = \rho_{i+1}^{(\beta)} \, \rho_{i}^{(\alpha)} \, \rho_{i+1}^{(\alpha)}$, for $0<\alpha<\beta$,
we obtain relation:
$$x_{ij}^{(\alpha)}x_{ik}^{(\alpha)}x_{jk}^{(\beta)}=x_{jk}^{(\beta)}x_{ik}^{(\alpha)}x_{ij}^{(\alpha)}.$$

\begin{proposition}
$\widetilde{M_kVH}_n$ is the quotient of $M_kVH_n$ by the relations
$$
x_{ij}^{(0)}x_{ik}^{(\alpha)}x_{jk}^{(\alpha)}=x_{ik}^{(\alpha)}x_{jk}^{(\alpha)}x_{ij}^{(0)},~~\alpha\neq 0,
$$
$$
x_{ij}^{(\alpha)}x_{ik}^{(\alpha)}x_{jk}^{(\beta)}=x_{jk}^{(\beta)}x_{ik}^{(\alpha)}x_{ij}^{(\alpha)},~~0<\alpha<\beta.
$$
\end{proposition}


\section{Some quotients and subgroups of multi-virtual braid group}
\label{QG}
\medskip

If we add to $M_kVB_n$ the forbidden relation $F1$ for all $i = 1, 2, \ldots, n-2,~\alpha = 0, 1, \ldots, k-1,$ we get the {\it multi-welded braid group} $M_kWB_n$. 
If we add to $M_kVB_n$ the forbidden relation $F1$ and $F2$ for all $i = 1, 2, \ldots, n-2,~\alpha = 0, 1, \ldots, k-1,$ we get the {\it multi-unrestricted  braid group} $M_kUB_n$. 

In the previous section we defined  two maps:
$$
\varphi_{n,k} \colon M_kVB_n \to S_{n},~~\sigma_i \mapsto \rho_i,~~\rho_i^{(\alpha)} \mapsto \rho_i,~~i = 1, 2, \ldots, n-2,~\alpha = 0,1, 2, \ldots, k-1,  
$$ 
with the kernel $M_kVP_n = \ker \varphi_{n,k}$, the multi-virtual pure  braid group, 
and the  map
$$
\psi_{n,k} \colon M_kVB_n \to S_{n},~~\sigma_i \mapsto e,~~\rho_i^{(\alpha)} \mapsto \rho_i,~~i = 1, 2, \ldots, n-2,~\alpha = 0,1, 2, \ldots, k-1,  
$$ 
with the kernel $M_kVH_n = \ker \psi_{n,k}$, the multi-virtual  semi-pure braid group.

The first  map induces  the map
$$
\varphi_{n,k} \colon M_kWB_n \to S_{n},~~\sigma_i \mapsto \rho_i,~~\rho_i^{(\alpha)} \mapsto \rho_i,~~i = 1, 2, \ldots, n-2,~\alpha = 0, 1, 2, \ldots, k-1.  
$$ 
The kernel $\ker \varphi_{n,k}$ is the {\it  multi-welded pure  braid group} $M_kWP_n$.
The second map induces  the map
$$
\psi_{n,k} \colon M_kWB_n \to S_{n},~~\sigma_i \mapsto e,~~\rho_i^{(\alpha)} \mapsto \rho_i,~~i = 1, 2, \ldots, n-2,~\alpha = 0, 1, 2, \ldots, k-1.  
$$ 
The kernel $\ker \psi_{n,k}$ is the {\it  multi-welded semi-pure braid group} $M_kWH_n$.

Also, the first  map induces the   map
$$
\varphi_{n,k} \colon M_kUB_n \to S_{n},~~\sigma_i \mapsto \rho_i,~~\rho_i^{(\alpha)} \mapsto \rho_i,~~i = 1, 2, \ldots, n-2,~\alpha = 0, 1, 2, \ldots, k-1.  
$$ 
The kernel $\ker \varphi_{n,k}$ is the {\it  multi-unrestricted pure  braid group} $M_kUP_n$.
The second map induces the map
$$
\psi_{n,k} \colon M_kUB_n \to S_{n},~~\sigma_i \mapsto e,~~\rho_i^{(\alpha)} \mapsto \rho_i,~~i = 1, 2, \ldots, n-2,~\alpha = 0, 1, 2, \ldots, k-1.  
$$ 
The kernel $\ker \psi_{n,k}$ is the {\it  multi-unrestricted semi-pure braid group} $M_kUH_n$.

The generators and defining relations for the pure braid groups $M_kWP_n$ and $M_kUP_n$, and for semi-pure $M_kWH_n$ and $M_kUH_n$ are stated as follows.
The group $M_kWB_n$ is quotient of the group $M_kVB_n$ by relation $F1$.

From  the relation $\sigma_i\sigma_{i+1}\rho_{i}^{(\alpha)}=\rho_{i+1}^{(\alpha)}\sigma_i\sigma_{i+1}$, we obtain the following relation:
$$\lambda_{ij}^{(0)}\lambda_{ik}^{(0)}\lambda_{jk}^{(\alpha)}=\lambda_{jk}^{(\alpha)}\lambda_{ik}^{(0)}\lambda_{ij}^{(0)} \text{, for } \alpha\neq 0.$$
 $$\lambda_{ij}^{(0)}\lambda_{ik}^{(0)}=\lambda_{ik}^{(0)}\lambda_{ij}^{(0)}.$$

The group $M_kUP_n$ is the quotient of group $M_kWP_n$ by relation $F2$. From  the relation $\sigma_{i+1}\sigma_{i}\rho_{i+1}^{(\alpha)}=\rho_{i}^{(\alpha)}\sigma_{i+1}\sigma_{i}$, we obtain the following relation:
$$\lambda_{jk}^{(0)}\lambda_{ik}^{(0)}\lambda_{ij}^{(\alpha)}=\lambda_{ij}^{(\alpha)}\lambda_{ik}^{(0)}\lambda_{jk}^{(0)} \text{, for } \alpha\neq 0,$$
 $$\lambda_{jk}^{(0)}\lambda_{ik}^{(0)}=\lambda_{ik}^{(0)}\lambda_{jk}^{(0)}.$$
Similarly, we can find presentation for the semi-pure braid groups $M_kWH_n$ and $M_kUH_n$.\\

The group $M_kWB_n$ is the quotient of the group $M_kVB_n$ by relation $F1$.\\
From the relation $\sigma_i\sigma_{i+1}\rho_{i}^{(\alpha)}=\rho_{i+1}^{(\alpha)}\sigma_i\sigma_{i+1}$, we obtain the following relations:

$$x_{ik}^{(0)}x_{ji}^{(0)}x_{ki}^{(\alpha)}=x_{ij}^{(\alpha)}x_{jk}^{(0)}x_{ij}^{(0)} \text{, for } \alpha\neq 0,$$
 $$x_{jk}^{(0)}x_{ik}^{(0)}=x_{ik}^{(0)}x_{jk}^{(0)}.$$
The group $M_kUH_n$ is the quotient of group $M_kWH_n$ by relation $F2$. From  the relation $\sigma_{i+1}\sigma_{i}\rho_{i+1}^{(\alpha)}=\rho_{i}^{(\alpha)}\sigma_{i+1}\sigma_{i}$, we obtain the following relations:\\
$$x_{ik}^{(0)}x_{kj}^{(0)}x_{ki}^{(\alpha)}=x_{jk}^{(\alpha)}x_{ij}^{(0)}x_{jk}^{(0)} \text{, for } \alpha\neq 0,$$
 $$x_{jk}^{(0)}x_{ik}^{(0)}=x_{ik}^{(0)}x_{jk}^{(0)}.$$



\section*{Concluding Remarks}  \label{CON}

This work opens several directions for future research within the framework of multi-virtual braid theory. Below are some open problems and potential avenues for further investigation:

\begin{enumerate}
    \item Construct representations of the multi-virtual braid group using automorphisms of groups and quandles.

    \item Investigate which algebraic structures correspond to braid-like groups. For instance, it is well known that the Hecke algebra corresponds to the classical braid group $B_n$. What analogous structures, if any, correspond to  multi-virtual braid groups?

    \item Define a group and a quandle for multi-virtual knots and establish that these structures are indeed invariants under appropriate equivalence moves.
    
    \item Loday and Stein \cite{LS} introduced parametric braid groups. Is it possible generalized this construction to define multi-virtual braid groups?

\end{enumerate}

\section*{Acknowledgments}
The first and second authors are supported by the Ministry of Science and Higher Education of Russia (agreement No. 075-02-2026-1339). The third author gratefully acknowledges NBHM for the postdoctoral fellowship (Ref. No.: 0204/21(1)/2025-R\& D-II/16279). The fourth author would like to acknowledge the support provided by the project CRG/2023/004921.  

\section*{Conflict of Interest}
The authors declare no conflicts of interest.

\end{document}